\documentclass{article}
\usepackage{amsfonts,graphicx,lscape}

\usepackage{amssymb}
\usepackage{eucal}
\usepackage{amsmath}
\usepackage{amsthm}
\usepackage{amscd}
\usepackage{graphics}
\usepackage{graphicx}
\usepackage{amsfonts}
\usepackage{latexsym}
\usepackage[all]{xy}

\numberwithin{equation}{section}
\newtheorem{thm}{\bf Theorem}[section]

\newtheorem{prop}[thm]{\bf Proposition}

\theoremstyle{remark}
\newtheorem{rem}{\bf Remark}[section]

\newcommand{\bc}{\begin{center}}
\newcommand{\ec}{\end{center}}
\newcommand{\bec}{\begin{equation}}
\newcommand{\eec}{\end{equation}}
\newcommand{\bea}{\begin{eqnarray}}
\newcommand{\eea}{\end{eqnarray}}
\newcommand{\ba}{\begin{array}}
\newcommand{\ea}{\end{array}}

\def\ds{\displaystyle}
\def\vs{\vspace{4pt}}

\begin{document}

\title{Symplectic realizations and symmetries of a Lotka-Volterra type system}
\author{Cristian L\u azureanu and Tudor B\^inzar\\
{\small Department of Mathematics, "Politehnica" University of Timi\c soara}\\
{\small Pia\c ta Victoriei nr. 2, 300006 Timi\c soara, Rom\^ania}\\
{\small E-mail: cristian.lazureanu@upt.ro;tudor.binzar@upt.ro}}
\date{}

\maketitle

\begin{abstract}
\noindent In this paper a Lotka-Volterra type system is considered. For such
a system, bi-Hamiltonian formulation, symplectic realizations
 and symmetries are presented.
\footnote{REGULAR \& CHAOTIC DYNAMICS,  Vol. 18,  Issue: 3 (2013)  Pages: 203--213,   DOI: 10.1134/S1560354713030015   }
\end{abstract}


\section{Introduction}
The dynamical systems of Lotka-Volterra type have significant importance in biology for interaction
models of biological species \cite{Lotka}, \cite{Volterra}, in chemistry for autocatalytic
chemical reactions \cite{Chenciner}, \cite{Toth}, in hydrodynamics \cite{Busse}, \cite{GiReZa}, and so on.

These systems have been widely investigated from different points
of view. We mention some studied topics, namely: integrals and
invariant manifold \cite{CaiFei}, \cite{Labrunie}, \cite{LeaMir},
\cite{LliVal},  Hamiltonian structure \cite{HerFai}, \cite{Plank},
\cite{TudGir}, \cite{BalBlaMus}, symmetries \cite{AlMaMo},
\cite{GonPer}, stability \cite{LeaMir}, \cite{May}, numerical
integration \cite{PopAro}, Kowalevski–-Painleve property
\cite{ConDam}, and many others.

In our paper, the following Lotka-Volterra type system
\begin{equation}\label{eq1.1}
\left\{\begin{array}{l}
\dot x=x(by+cz)\\
\dot y=y(ax-by+cz+d)\\
\dot z=z(-ax+by-cz-d)
\end{array}\right.,
\end{equation}
where $a,c\in{\bf R}^*,b,d\in{\bf R}$, is considered and some symmetries are given.

For our purposes, a Hamilton-Poisson realization and a symplectic realization of system (\ref{eq1.1}) are required.

Theoretical details about symmetries of differential equations can be found in \cite{BluKum}, \cite{Dam}
, \cite{FokFuc}, \cite{Fuc}, \cite{Olver}.

A similar study for Maxwell-Bloch equations was presented by P.A.Damianou and P.G.Paschali in \cite{DamPas}.

\section{Hamiltonian structures and symmetries for considered system in the case $d=0$}
In this section, we consider system (\ref{eq1.1}) with $d=0$, i.e.
\begin{equation}\label{eq2.1}
\left\{\begin{array}{l}
\dot x=x(by+cz)\\
\dot y=y(ax-by+cz)\\
\dot z=z(-ax+by-cz)
\end{array}\right.
\end{equation}

A bi-Hamiltonian structure, a symplectic realization and some symmetries of system (\ref{eq2.1}) are given.

For system (\ref{eq2.1}), the functions $H_1,H_2\in{\cal C}^\infty
({\bf R}^3,{\bf R}),$
$$~H_1(x,y,z)=yz$$
and
$$H_2(x,y,z)=x(ax-2by+2cz)$$
are constants of motion.

Let us consider the linear Poisson algebra ${\cal P}_1$,
\begin{eqnarray}
\{x,y\}_1&=&\alpha_1x+\alpha_2y+\alpha_3z,\nonumber\\
\{x,z\}_1&=&\beta_1x+\beta_2y+\beta_3z,\nonumber\\
\{y,z\}_1&=&\gamma_1x+\gamma_2y+\gamma_3z.\nonumber
\end{eqnarray}

Imposing the condition that $C=H_1$ to be a Casimir for ${\cal P}_1$, it results $\alpha_1=\alpha_3=\beta_1=\beta_2=\gamma_1=\gamma_2=\gamma_3=0$,
$\beta_3=-\alpha_2.$ If the Hamiltonian function is $H=H_2$, we get the following dynamical system:
\begin{equation}\label{eq2.02}
\left\{\begin{array}{l}
\dot x=-2\alpha_2x(by+cz)\\
\dot y=-2\alpha_2y(ax-by+cz)\\
\dot z=-2\alpha_2z(-ax+by-cz)
\end{array}\right.
\end{equation}

Taking $\alpha_2=-\frac{1}{2}$, the above system is the considered system (\ref{eq2.1}). Thus,
$$\{x,y\}_1=-\frac{1}{2}y,~\{x,z\}_1=\frac{1}{2}z,~\{y,z\}_1=0,$$
or in coordinates, using matrix notation,
$$\pi_1(x,y,z)=\left[\begin{array}{ccc}
0&-\ds\frac{1}{2}y&\ds\frac{1}{2}z\vs\\
\ds\frac{1}{2}y&0&0\vs\\
-\ds\frac{1}{2}z&0&0\end{array}\right].$$
Therefore we consider the three-dimensional Lie algebra $g_1$ given by
$$[E_1,E_2]=-\ds\frac{1}{2}E_2,[E_1,E_3]=\ds\frac{1}{2}E_3,[E_2,E_3]=0,$$
where
$$E_1=\left[\begin{array}{rrr}
-\ds{\frac{1}{2}}&0&0\\
0&\ds{\frac{1}{2}}&0\\
0&0&0\end{array}\right],E_2=\left[\begin{array}{rrr}
0&0&1\\0&0&0\\
0&0&0\end{array}\right],E_3=\left[\begin{array}{rrr}
0&0&0\\0&0&1\\
0&0&0\end{array}\right].$$

As a real vector space, $g_1$ is generated by the base $B_{g_1}=\{E_1,E_2,E_3\},$ whence
$$g_1=\{X\in gl(3,{\bf R})|~X=\left[\begin{array}{ccc}
-\alpha &0&\beta \\
0&\alpha &\gamma\\
0&0&0
\end{array}\right],~\alpha ,\beta ,\gamma \in{\bf R}\}.$$

In order to provide the Lie group $G_1$ generated by the Lie algebra $g_1$, we consider
$$A\!=\!\exp (k E_1)\cdot\exp (l E_2)\cdot\exp(m E_3)\!=\!\left[\!\begin{array}{ccc}
e^{-\frac{1}{2}k}&0&l e^{-\frac{1}{2}k }\\
0&e^{\frac{1}{2}k}&m e^{\frac{1}{2}k}\\
0&0&1
\end{array}\!\right]$$
Taking $k=2u$, $l=ve^u$, $m=we^{-u}$, it follows
$$G_1=\{A\in GL(3,{\bf R})|~A=\left[\begin{array}{ccc}
e^{-u}&0&v\\
0&e^u&w\\
0&0&1
\end{array}\right],~u,v,w\in{\bf R}\}.$$

In the same manner, interchanging $H_1$ to $H_2$, we get the linear Poisson algebra ${\cal P}_2$:
$$\{x,y\}_2=cx,~\{x,z\}_2=bx,~\{y,z\}_2=ax-by+cz,$$
or
$$\pi_2(x,y,z)=\left[\begin{array}{ccc}0&cx&bx\\
-cx&0&ax-by+cz\\
-bx&-ax+by-cz&0\end{array}\right].
$$
Let us consider the 3D Lie algebra $g_2$, given by
\begin{eqnarray}
&&[X_1,X_2]=cX_1,[X_1,X_3]=bX_1,\nonumber\\
&&[X_2,X_3]=aX_1-bX_2+cX_3,\nonumber \end{eqnarray}
where
$$X_1=\left[\begin{array}{rrr}
0&1&0\\
0&0&0\\
0&0&0\end{array}\right]\!,X_2=\left[\begin{array}{rrr}
-c&\ds\frac{a}{b}&c\vs\\
0&0&0\\
0&\ds\frac{a}{b}&c\end{array}\right]\!,X_3=\left[\begin{array}{rrr}
-b&0&b\\
0&0&0\\
0&0&b\end{array}\right]\!,$$
$b\not=0.$

The Lie group $G_2$ generated by the Lie algebra $g_2$ has the generic element
\begin{eqnarray}
A&=&\exp (\alpha X_1)\cdot\exp (\beta X_2)\cdot\exp(\gamma X_3)=\nonumber\end{eqnarray}
$$=\!\left[\!\begin{array}{ccc}
e^{-\beta c-\gamma b}&\alpha+\frac{a}{2bc}(e^{\beta c}-e^{-\beta c})&\frac{1}{2}(e^{\beta c+\gamma b}-e^{-\beta c-\gamma b})\\
0&1&0\\
0&\frac{a}{bc}(e^{\beta c}-1)&e^{\beta c+\gamma b}
\end{array}\!\right]$$

Taking $\alpha =u, e^{\beta }=v,e^{\gamma }=w$, we obtain
$$G_2=\{~
A=\left[\begin{array}{ccc}
v^{-c}w^{-b}&u+\frac{a}{2bc}(v^c-v^{-c})&\frac{1}{2}(v^cw^b-v^{-c}w^{-b})\\
0&1&0\\
0&\frac{a}{bc}(v^c-1)&v^cw^b
\end{array}\right],~u,v,w\in{\bf R},v>0,w>0\}.$$

In the case $b=0$, the base $\{Y_1,Y_2,Y_3\}$ of $g_2$ is given by
$$Y_1=\left[\begin{array}{rrr}
0&1&0\\
0&0&0\\
0&0&0\end{array}\right],Y_2=\left[\begin{array}{rrr}
c&a&-c\vs\\
0&2c&0\\
0&a&0\end{array}\right],Y_3=\left[\begin{array}{rrr}
0&0&-1\\
0&0&0\\
0&0&0\end{array}\right]$$
and the corresponding Lie group is
$$G_2^{b=0}=\{~A=\left[\begin{array}{ccc}
e^{\beta c}&\frac{a}{2c}(e^{2\beta c}-1)+\alpha e^{2\beta c}&1-(1+\gamma)e^{\beta c}\\
0&e^{2\beta c}&0\\0&\frac{a}{2c}(e^{2\beta c}-1)&1
\end{array}\right],~\alpha ,\beta , \gamma\in{\bf R}\}.$$

Since
$$\pi_1\cdot\nabla H_2=\pi_2\cdot\nabla H_1=\left(\begin{array}{c}x(by+cz)\\y(ax-by+cz)\\z(-ax+by-cz)\end{array}\right)=
\left(\begin{array}{c}\dot{x}\\\dot{y}\\\dot{z}\end{array}\right),$$
system (\ref{eq2.1}) is a bi-Hamiltonian system. For the $\pi_1$ bracket, $H_2$ is the Hamiltonian and
$H_1$ is a Casimir. For the $\pi_2$ bracket, $H_1$ is the Hamiltonian and $H_2$ is a Casimir.

We recall that for a system $\dot x=f(x)$, where $f:M\to TM$, and $M$ is a smooth manifold of finite dimension, a vector
field ${\bf X}$ is called:

$\bullet $ {\em a symmetry} if $\ds{\frac{\partial {\bf X}}{\partial t}+[{\bf X},{\bf X}_f]=0}$, where ${\bf X}_f$ is
the vector field defined by the system;

$\bullet $ {\em a Lie-point symmetry} if its first prolongation transforms solutions of the system into other solutions;

$\bullet $ {\em a conformal symmetry} if the Lie derivative along ${\bf X}$ satisfies $L_{\bf X}\pi =\lambda\pi $ and
$L_{\bf X}H=\nu H$, for some scalars $\lambda ,\nu $, where the Poisson tensor $\pi $ and the Hamiltonian $H$ give the
Hamilton-Poisson realization of the system;

$\bullet $ {\em a master symmetry} if $[[{\bf X},{\bf X}_f],{\bf X}_f]=0$, but $[{\bf X},{\bf X}_f]\not=0$.

First, we give a characterization of the vector field
$$~\ds{{\bf X}=\alpha t\frac{\partial }{\partial t}+x\frac{\partial }{\partial x}+
y\frac{\partial }{\partial y}+z\frac{\partial }{\partial z}},~\alpha\in{\bf R},~$$
to be a Lie-point symmetry for system
\begin{equation}\label{eq2.2}
\left\{\begin{array}{l}
\dot x=f(x,y,z)\\
\dot y=g(x,y,z)\\
\dot z=h(x,y,z)
\end{array}\right.
\end{equation}

\begin{thm}\label{th2.1}
Let $f,g,h$ be three real functions of class $C^1$ on a cone in ${\bf R}^3$.

The vector field $$~\ds{{\bf X}=\alpha t\frac{\partial }{\partial t}+x\frac{\partial }{\partial x}+
y\frac{\partial }{\partial y}+z\frac{\partial }{\partial z}},~\alpha\in{\bf R},~$$ is a Lie-point symmetry of system (\ref{eq2.2})
if and only if $f,g,h$ are homogeneous functions of degree $1-\alpha $.
\end{thm}

{\bf Proof.} The vector field ${\bf X} $ is a Lie point symmetry
for system (\ref{eq2.2}) if and only if its first prolongation
$$pr^{(1)}({\bf X})={\bf X}+(\dot{x}-\alpha\dot x)\frac{\partial }{\partial \dot x}+
(\dot{y}-\alpha\dot y)\frac{\partial }{\partial \dot y}+(\dot{z}-\alpha\dot z)\frac{\partial }{\partial\dot z}$$
applied to the system equations vanishes. This condition is equivalent to
$$\left\{\begin{array}{l}
\ds{(1-\alpha)f(x,y,z)-x\frac{\partial f}{\partial x}-y\frac{\partial f}{\partial y}-z\frac{\partial f}{\partial z}=0}\vs\\
\ds{(1-\alpha)g(x,y,z)-x\frac{\partial g}{\partial x}-y\frac{\partial g}{\partial y}-z\frac{\partial g}{\partial z}=0}~~,\vs\\
\ds{(1-\alpha)h(x,y,z)-x\frac{\partial h}{\partial
x}-y\frac{\partial h}{\partial y}-z\frac{\partial h}{\partial
z}=0}\end{array}\right.$$ that is $f,g,h$ are homogeneous
functions of degree $1-\alpha $.\vs

 The next result one furnishes a Lie point symmetry of system (\ref{eq2.1}) and a
 conformal symmetry.
\begin{prop} The vector field
$$~\ds{{\bf X}=-t\frac{\partial }{\partial t}+x\frac{\partial }{\partial x}+
y\frac{\partial }{\partial y}+z\frac{\partial }{\partial z}}~$$
is a Lie point symmetry of system (\ref{eq2.1}). Moreover, ${\bf X}$ is a conformal symmetry.
\end{prop}

{\bf Proof.}  Using {\bf Theorem} \ref{th2.1}, it follows that
${\bf X}$ is a Lie point symmetry of system (\ref{eq2.1}).

One can easily check that
$$L_{\bf X}\pi_1=-\pi_1,L_{\bf X}\pi_2=-\pi_2,L_{\bf X}H_1=2H_1,
L_{\bf X}H_2=2H_2,$$ whence ${\bf X}$ is a conformal symmetry. \vs

The following result provides a master symmetry of our considered system.

\begin{prop} The vector field
$$\begin{array}{lll}
\overrightarrow{X}&=&(k_1x+k_2bxy+k_2cxz)\ds\frac{\partial }{\partial x}+\vs\\
&&+(k_1y+k_2axy-k_2by^2+k_2cyz)\ds\frac{\partial }{\partial y}+\vs\\
&&+(k_1z-k_2axz+k_2byz-k_2cz^2)\ds\frac{\partial }{\partial z},\end{array}$$
where $k_1\in{\bf R}^*,~k_2\in{\bf R},$
is a master symmetry of system (\ref{eq2.1}).
\end{prop}

{\bf Proof.}  We denote by $\overrightarrow{V}$ the associated
vector field of system (\ref{eq2.1}), that is\\
$\ds{\overrightarrow{V}=(bxy+cxz)\frac{\partial }{\partial
x}+(axy-by^2+cyz)\frac{\partial }{\partial y}}+$\\
$+\ds{(-axz+byz-cz^2)\frac{\partial }{\partial z}}.$

It follows that the following relations
$$[\overrightarrow{X},\overrightarrow{V}]=k_1\overrightarrow{V},~\left[[\overrightarrow{X},\overrightarrow{V}],\overrightarrow{V}\right]=
\overrightarrow{0}$$
hold.

Therefore $\overrightarrow{X}$ is a master symmetry of system
(\ref{eq2.1}). \vs

In the following a symplectic realization of system (\ref{eq2.1}) is given. Using this fact, the symmetries of Newton's equations are studied.

The next theorem states that the system (\ref{eq2.1}) can be regarded as a Hamiltonian mechanical system.

\begin{thm} The Hamilton-Poisson mechanical system $({\bf R}^3,\pi_1,H_2)$ has a full symplectic realization
$({\bf R}^4,\omega ,\tilde{H})$, where $\omega =\mbox{d}p_1\wedge\mbox{d}q_1+\mbox{d}p_2\wedge\mbox{d}q_2$ and
$$\tilde{H}=\frac{1}{a}\left[p_1^2-\left(bp_2e^{\frac{a}{2}q_1}-ce^{-\frac{a}{2}q_1}\right)^2\right].$$
\end{thm}

{\bf Proof.} The corresponding Hamilton's equations are
\begin{equation}\label{eq2.3}
\left\{\begin{array}{l}
\dot{q}_1=\ds\frac{2}{a}p_1\vs\\
\dot{q}_2=-\ds\frac{2b^2}{a}p_2e^{aq_1}+\ds\frac{2bc}{a}\vs\\
\dot{p}_1=-\ds c^2e^{-aq_1}+b^2p_2^2e^{aq_1}\vs\\
\dot{p}_2=0
\end{array}\right.
\end{equation}

We define the application $\varphi :{\bf R}^4\to{\bf R}^3$ by $$
\varphi (q_1,q_2,p_1,p_2)=(x,y,z),$$ where
\begin{eqnarray}
x&=&\ds{\frac{1}{a}p_1+\frac{b}{a}}p_2e^{\frac{a}{2}q_1}-\ds\frac{c}{a}e^{-\frac{a}{2}q_1}\nonumber\\
y&=&p_2e^{\frac{a}{2}q_1}\nonumber\\
z&=&e^{-\frac{a}{2}q_1}.\nonumber\end{eqnarray}

It follows that $\varphi $ is a surjective submersion, the equations (\ref{eq2.3}) are mapped onto the equations (\ref{eq2.1}), the
canonical structure $\{.,.\}_{\omega }$ is mapped onto the Poisson structure $\pi_1$, as required.

We also remark that $H_2=\tilde{H}$ and $H_1=p_2.$ \vs

It is natural to ask which non-canonical bracket $\tilde{\omega }$ in ${\bf R}^4$ is mapped onto the $\pi_2$ bracket.

Considering now the Hamiltonian $H=p_2$ and taking into account the Jacobi's identity, it results
\begin{eqnarray}
\{q_1,p_2\}_{\tilde{\omega }}&=&\frac{2}{a}p_1,\nonumber\\
\{q_2,p_2\}_{\tilde{\omega }}&=&\frac{2b}{a}(c-bp_2e^{aq_1}),\nonumber\\
\{p_1,p_2\}_{\tilde{\omega }}&=&b^2p_2^2e^{aq_1}-c^2e^{-aq_1},\nonumber\\
\{q_1,p_1\}_{\tilde{\omega }}&=&\frac{2b}{a}(bp_2e^{aq_1}-c),\nonumber\\
\{q_1,q_2\}_{\tilde{\omega }}&=&F(p_1,p_2,q_1,q_2),\nonumber\\
\{p_1,q_2\}_{\tilde{\omega }}&=&G(p_1,p_2,q_1,q_2),\nonumber
\end{eqnarray}
where the functions $F$ and $G$ satisfies the relations:
$$-\frac{2}{a}p_1\frac{\partial F}{\partial q_1}+\frac{2b}{a}(bp_2e^{aq_1}-c)\frac{\partial F}{\partial q_2}+
(c^2e^{-aq_1}-b^2p_2^2e^{aq_1})\frac{\partial F}{\partial p_1}+\frac{2}{a}G=\frac{4b^2}{a^2}p_1e^{aq_1}$$
$$\frac{2}{a}p_1\frac{\partial G}{\partial q_1}-\frac{2b}{a}(bp_2e^{aq_1}-c)\frac{\partial G}{\partial q_2}-
(c^2e^{-aq_1}-b^2p_2^2e^{aq_1})\frac{\partial G}{\partial p_1}-a(b^2p_2^2e^{aq_1}+c^2e^{-aq_1})F=~~~~~~~~~~~~$$
$$~~~~~~~~~~~~~~~~~~~~~~~~~~~~~~~~~~=\frac{6b^4}{a}p_2^2e^{2aq_1}-\frac{8b^3c}{a}p_2e^{aq_1}+\frac{2b^2c^2}{a}$$
$$\frac{2b(c-bp_2e^{aq_1})}{a}\left(\frac{\partial F}{\partial q_1}+\frac{\partial G}{\partial p_1}\right)+G\frac{\partial F}{\partial q_2}-
(c^2e^{-aq_1}-b^2p_2^2e^{aq_1})\frac{\partial F}{\partial p_2}-F\frac{\partial G}{\partial q_2}-
\frac{2p_1}{a}\frac{\partial G}{\partial p_2}+
2b^2p_2e^{aq_1}F=~~~~~~~~~~~~$$$$~~~~~~~~~~~~~~~~~~~~~~~=\frac{4b^3}{a^2}(ce^{aq_1}-bp_2e^{2aq_1}).$$

In the particular case $b=0$ one finds $$F_1(p_1,p_2,q_1,q_2)=p_1~,~G_1(p_1,p_2,q_1,q_2)=-\frac{ac^2}{2}e^{-aq_1}$$
and respectively
$$F_2(p_1,p_2,q_1,q_2)=p_1+e^{-aq_1}~,$$$$G_2(p_1,p_2,q_1,q_2)=-\frac{a}{2}e^{-aq_1}(2p_1+c^2).$$
Thus, we obtain two non-canonical brackets $\tilde{\omega}_1$ and $\tilde{\omega}_2$.

The first one is compatible with $\omega $ but it is degenerate.

A symmetry of system (\ref{eq2.3}) is
$$Z_0=-t\frac{\partial }{\partial t}+p_1\frac{\partial }{\partial p_1}+2p_2\frac{\partial }{\partial p_2}
-\frac{2}{a}\frac{\partial }{\partial q_1}+(\frac{3a}{2}p_2-q_2)\frac{\partial }{\partial q_2}$$
which satisfies
$$L_{Z_0}(J_0)=-J_0~,~L_{Z_0}(J_1)=-J_1~,~L_{Z_0}(\tilde{H})=2\tilde{H},$$
where $J_0$ and $J_1$ are respectively Poisson tensors associated with $\omega $ and $\tilde{\omega }_1$.

Defining the operator ${\cal R}=J_1J_0^{-1}$, we get a symmetry
$$Z_1={\cal R}Z_0=\frac{1}{2}c^2e^{-aq_1}(2q_2-ap_2)\frac{\partial }{\partial p_1}+\frac{2}{a}(p_1^2-c^2e^{-aq_1})\frac{\partial }{\partial p_2}+
\frac{p_1}{a}(ap_2-2q_2)\frac{\partial }{\partial q_1}+(p_1^2-c^2e^{-aq_1})\frac{\partial }{\partial q_2}$$
of system (\ref{eq2.3}). This symmetry leads to a symmetry of system (\ref{eq2.1}) for $b=0$:
$$\overrightarrow{Z}=\frac{c}{4}xz(2k-ayz)\frac{\partial }{\partial x}-
\frac{1}{4}(4a\frac{x^2}{z}+8cx+a^2xy^2z+acy^2z^2-2akxy-2ckyz)\frac{\partial }{\partial y}
-\frac{z}{4}(ax+cz)(2k-ayz)\frac{\partial }{\partial z},$$
$k\in{\bf R},$ 
which sending $\pi_1$ bracket to $\pi_2$ and $H_1$ to $H_2$. It follows that
$~\ds{\overrightarrow{Y}=\overrightarrow{Z}+\left(x\frac{\partial }{\partial x}+y\frac{\partial }{\partial y}+z\frac{\partial }{\partial z}
\right)}~$
is a master symmetry of system (\ref{eq2.1}) which sending $\pi_1$ bracket to $\pi_2-\pi_1$ and $H_1$ to $2H_1-H_2$.

The second bracket $\tilde{\omega }_2$ is non-degenerate but it is not compatible with $\omega $ and they do not generate a
recursion operator.\vs

In the sequel we find the symmetries of Newton's equations. In the case $b\not=0$, from Hamilton's equations (\ref{eq2.3})
one obtains Newton's equations:
\begin{equation}\label{eq2.4}
\ddot{q}_1-\ds\frac{a}{2b^2}e^{-aq_1}\dot{q}_2^2+\ds\frac{2c}{b}e^{-aq_1}\dot{q}_2=0
\end{equation}
\begin{equation}\label{eq2.5}
\ddot{q}_2-a\dot{q}_1\dot{q}_2+2bc\dot{q}_1=0.
\end{equation}

These are also Lagrange's equations generated by the Lagrangian
$$L=\frac{a}{4}\dot{q}_1^2+\left(\ds\frac{c}{b}\dot{q}_2-\ds\frac{a}{4b^2}\dot{q}_2^2\right)e^{-aq_1}.$$

A vector field $$\overrightarrow{v}=\xi (q_1,q_2,t)\frac{\partial }{\partial t}+\eta_1(q_1,q_2,t)\frac{\partial }{\partial q_1}+
\eta_2(q_1,q_2,t)\frac{\partial }{\partial q_2}$$
is a Lie-point symmetry for Newton's equations if the action of its second prolongation on Newton's equations vanishes. Thus, for the equation
(\ref{eq2.5}), the following condition is obtained:
$$\begin{array}{l}
(-2bc\eta_{1,t}-\eta_{2,tt})+\dot{q}_1(-2bc\eta_{1,q_1}-2bc\xi_t+a\eta_{2,t}-2\eta_{2,tq_1}+2bc\eta_{2,q_2})+
\dot{q}_1^2(-2bc\xi_{q_1}+a\eta_{2,q_1}-\eta_{2,q_1q_1})+\vs\\
+\dot{q}_2\left(a\eta_{1,t}-2bc\eta_{1,q_2}-2\eta_{2,tq_2}+\ds\frac{2c}{b}e^{-aq_1}\eta_{2,q_1}+\xi_{tt}\right)+\vs\\
+\dot{q}_2^2\left(a\eta_{1,q_2}-\eta_{2,q_2q_2}-\ds\frac{a}{2b^2}e^{-aq_1}\eta_{2,q_1}+2\xi_{tq_2}-\ds\frac{2c}{b}e^{-aq_1}\xi_{q_1}\right)
+\dot{q}_2^3\left(\xi_{q_2q_2}+\ds\frac{a}{2b^2}e^{-aq_1}\xi_{q_1}\right)+\vs\\
+\dot{q}_1\dot{q}_2(a\eta_{1,q_1}-2\eta_{2,q_1q_2}+2\xi_{tq_1}-4bc\xi_{q_2})+
\dot{q}_1^2\dot{q}_2\xi_{q_1q_1}+\dot{q}_1\dot{q}_2^2(2\xi_{q_1q_2}+a\xi_{q_2})=0.\end{array}$$

The above equation must be satisfied identically in $t,$ $q_1,$
$q_2,$ $\dot{q}_1,$ $\dot{q}_2$, that are all independent. Doing
standard manipulation, we obtain:
$$\begin{array}{l}
\eta_{2}=\eta_2(q_2,t)~,~~
\eta_{1}=-\ds{\frac{1}{2bc}\cdot\frac{\partial \eta_2}{\partial t}+k_1}\vs\\
\xi =\xi(t),~\xi{}'(t)=\ds{\frac{a}{2bc}\cdot\frac{\partial \eta_2}{\partial t}+\frac{\partial \eta_2}{\partial q_2}},
\end{array}$$
where $k_1$ is an arbitrary real constant.

Taking into account the above result, the second prolongation of $\overrightarrow{v}$ on equation (\ref{eq2.4}) gives us:
\begin{eqnarray}
&&\dot{q}_2^2\left(-\ds\frac{a^2}{2b^2}e^{-aq_1}\eta_1+\ds\frac{a}{b^2}e^{-aq_1}\eta_{2,q_2}-\eta_{1,q_2q_2}\right)-
\dot{q}_1\dot{q}_2a\eta_{1,q_2}+\dot{q}_1(\xi_{tt}+2bc\eta_{1,q_2})+\vs\nonumber\\
&&+\dot{q}_2\left(\ds\frac{2ac}{b}e^{-aq_1}\eta_1+\ds\frac{a}{b^2}e^{-aq_1}\eta_{2,t}-\ds\frac{2c}{b}e^{-aq_1}\eta_{2,q_2}-
\ds\frac{2c}{b}e^{-aq_1}\xi_t-2\eta_{1,tq_2}\right)+\left(-\eta_{1,tt}-\ds\frac{2c}{b}e^{-aq_1}\eta_{2,t}\right)=0.\nonumber\end{eqnarray}

We get the overall result:
$$\left\{\begin{array}{l}\xi =a\alpha t+\beta\\\eta_1=2\alpha \\\eta_2=a\alpha q_2+\gamma\end{array}\right.~,$$
where $\alpha ,\beta ,\gamma $ are real constants.

Now, we conclude the following result:
\begin{thm}
The symmetries of Newton's equations are given by
$$\overrightarrow{v}=(\alpha at+\beta )\frac{\partial }{\partial t}+2\alpha\frac{\partial }{\partial q_1}+
(\alpha aq_2+\gamma)\frac{\partial }{\partial q_2},$$
where $\alpha,\beta,\gamma\in{\bf R}.$
\end{thm}

\begin{rem}
(\emph{i}) For $\alpha =\gamma =0$ and $\beta\not=0$, we have $\overrightarrow{v_1}=\beta\ds\frac{\partial }{\partial t}$ that represents the time translation symmetry which generates the conservation of energy $H$.

(\emph{ii}) For $\alpha =\beta =0$ and $\gamma\not=0$, we have $\overrightarrow{v_2}=\gamma\ds\frac{\partial }{\partial q_2}$ that represents a translation in the cyclic $q_2$ direction which is related to the conservation of $p_2$.

Moreover, using the Lagrangian $L$ and Noether's theory we deduce that both $\overrightarrow{v_1}$ and $\overrightarrow{v_2}$ are variational
symmetries since they satisfy the condition $~pr^{(1)}(\overrightarrow{v})L+L\mbox{div}(\xi )=0.$
\end{rem}

\begin{rem}
The 3-dimensional Lie algebra corresponding to the symmetries of Newton's equations endowed with the standard Lie bracket vector fields is
generated by the base $\{\overrightarrow{u}_1,\overrightarrow{u}_2,\overrightarrow{u}_3\}$, where
$$\begin{array}{l}
\overrightarrow{u}_1=\ds{-t\cdot\frac{\partial }{\partial t}-\frac{2}{a}\cdot\frac{\partial }{\partial q_1}-q_2\cdot\frac{\partial }{\partial q_2}}\vs\\
\ds{\overrightarrow{u}_2=\frac{\partial }{\partial t}}\vs\\
\ds{\overrightarrow{u}_3=\frac{\partial }{\partial q_2}}~.
\end{array}$$

The following relations
$$[\overrightarrow{u}_1,\overrightarrow{u}_2]=\overrightarrow{u}_2~,~~[\overrightarrow{u}_1,\overrightarrow{u}_3]=\overrightarrow{u}_3~,~~
[\overrightarrow{u}_2,\overrightarrow{u}_3]=\overrightarrow{0}$$
hold. Therefore this Lie algebra is of type V in Bianchi classification \cite{Bianchi}.
\end{rem}

\section{Hamiltonian structures and symmetries for considered system in the case $d\not=0$}
Let us consider system (\ref{eq1.1}) in the case $d\not=0$.

In this section a symplectic realization of system (\ref{eq1.1}) is given. Using this fact, the symmetries of Newton's equations are studied.

For our purpose we can use the same Hamilton-Poisson realization $({\bf R}^3,\pi_1,H_2)$ as like as in the case $d=0$, but for the sake of
simplicity we choose another realization. We consider the following Lie-Poisson structure:
$$\pi (x,y,z)=\left[\begin{array}{ccc}0&\ds\frac{c}{ad}(ax-2by)&\ds\frac{b}{ad}(ax+2cz)\vs\\
-\ds\frac{c}{ad}(ax-2by)&0&\ds\frac{1}{d}(ax-by+cz+d)\vs\\
-\ds\frac{b}{ad}(ax+2cz)&-\ds\frac{1}{d}(ax-by+cz+d)&0\end{array}\right].
$$
The Hamiltonian $H$ is given by
$$H(x,y,z)=dyz,$$ and
moreover, the function $C$,
$$C(x,y,z)=x(ax-2by+2cz+2d)-\frac{4bc}{a}yz,$$
is a Casimir of our configuration.

The next theorem states that the system (\ref{eq1.1}) can be regarded as a Hamiltonian mechanical system.

\begin{thm} The Hamilton-Poisson mechanical system $({\bf R}^3,\pi ,H)$ has a full symplectic realization
$$({\bf R}^4,\omega ,\tilde{H}),$$ where $$\omega =\mbox{d}p_1\wedge\mbox{d}q_1+\mbox{d}p_2\wedge\mbox{d}q_2$$ and
$$\tilde{H}=bcdq_1^2-\frac{1}{16bcd}\left(ap_2-p_1^2+4b^2c^2q_1^2\right)^2.$$
\end{thm}

{\bf Proof.} The corresponding Hamilton's equations are
\begin{equation}\label{eq3.1}
\left\{\begin{array}{l}
\dot{q}_1=\ds\frac{1}{4bcd}p_1\left(ap_2-p_1^2+4b^2c^2q_1^2\right)\vs\\
\dot{q}_2=-\ds\frac{a}{8bcd}\left(ap_2-p_1^2+4b^2c^2q_1^2\right)\vs\\
\dot{p}_1=-2bcdq_1+\ds\frac{bc}{d}\left(ap_2-p_1^2+4b^2c^2q_1^2\right)q_1\vs\\
\dot{p}_2=0
\end{array}\right.
\end{equation}

We define the application $\varphi :{\bf R}^4\to{\bf R}^3$ by
$$\varphi (q_1,q_2,p_1,p_2)=(x,y,z),$$
where
\begin{eqnarray}
x&=&\ds{\frac{1}{a}p_1+\frac{1}{2ad}\left(ap_2-p_1^2+4b^2c^2q_1^2\right)-\frac{d}{a}}\nonumber\\
y&=&cq_1+\ds\frac{1}{4bd}\left(ap_2-p_1^2+4b^2c^2q_1^2\right)\nonumber\\
z&=&bq_1-\ds\frac{1}{4cd}\left(ap_2-p_1^2+4b^2c^2q_1^2\right)\nonumber
\end{eqnarray}

It follows that $\varphi $ is a surjective submersion, the equations (\ref{eq3.1}) are mapped onto the equations (\ref{eq1.1}), the
canonical structure $\{.,.\}_{\omega }$ is mapped onto the Poisson structure $\pi $, as required.

We also remark that $H=\tilde{H}$ and $C=p_2.$ \vs

From Hamilton's equations (\ref{eq3.1}) we obtain by differentiation, Newton's equations:
$$4\ddot{q}_2\dot{q}_2-a^2q_1\dot{q}_1=0$$
$$4a^2\ddot{q}_1\dot{q}_2^2-64b^2c^2\dot{q}_2^4q_1-16abcd\dot{q}_2^3q_1-a^4\dot{q}_1^2q_1=0$$
These are also Lagrange's equations generated by the Lagrangian
$$L=\frac{a}{4}\frac{\dot{q}_1^2}{\dot{q}_2}+\frac{4bcd}{a^2}\dot{q}_2^2+\frac{4b^2c^2}{a}\dot{q}_2q_1^2+bcdq_1^2.$$

The condition for the vector field $$\overrightarrow{v}=\xi (q_1,q_2,t)\frac{\partial }{\partial t}+\eta_1(q_1,q_2,t)\frac{\partial }{\partial q_1}+\eta_2(q_1,q_2,t)\frac{\partial }{\partial q_2}$$
to be a Lie Point symmetry for Newton's equations leads to:
\noindent $4\dot{q}_2\ddot{\eta}_2-4\dot{q}_2^2\ddot{\xi }-a^2q_1\dot{\eta }_1+4\ddot{q}_2\dot{\eta }_2+
\dot{\xi}(a^2q_1\dot{q}_1-12\ddot{q}_2\dot{q}_2)-a^2\dot{q}_1\eta_1=0$\vs\\
$4a^2\dot{q}_2^2\ddot{\eta}_1-4a^2\dot{q}_1\dot{q}_2^2\ddot{\xi }-2a^4\dot{q}_1q_1\dot{\eta }_1+
(8a^2\ddot{q}_1\dot{q}_2-256b^2c^2\dot{q}_2^3q_1-
48abcd\dot{q}_2^2q_1)\dot{\eta }_2+(2a^4\dot{q}_1^2q_1-16a^2\ddot{q}_1\dot{q}_2^2+256b^2c^2\dot{q}_2^4q_1+$\vs\\
$+48abcd\dot{q}_2^3q_1)\dot{\xi }-(64b^2c^2\dot{q}_2^4+16abcd\dot{q}_2^3+a^4\dot{q}_1^2)\eta_1=0.$

The resulting equations obtained by expanding $\dot{\xi
},\ddot{\xi },\dot{\eta }_1,\ddot{\eta }_1,\dot{\eta
}_2,\ddot{\eta }_2$ and replacing $\ddot{q}_1$ and $\ddot{q}_2$
must be satisfied identically in $t,$ $q_1,$ $q_2,$ $\dot{q}_1,$
$\dot{q}_2$, that are all independent. Doing standard
manipulation, we get the overall result:
$$\left\{\begin{array}{l}\xi =k_1\\\eta_1=0\\\eta_2=k_2\end{array}\right.$$
where $k_1,k_2$ are real constants.

For $k_2=0$ and $k_1\not=0$, we have $\overrightarrow{v_1}=k_1\ds\frac{\partial }{\partial t}$ that represents the time translation symmetry
which generates the conservation of energy $\tilde{H}$.

For $k_1=0$ and $k_2\not=0$, we have $\overrightarrow{v_2}=k_2\ds\frac{\partial }{\partial q_2}$ that represents a translation in the
cyclic $q_2$ direction which is related to the conservation of $p_2$.

Moreover, using the Lagrangian $L$ and Noether's theory we deduce that both $\overrightarrow{v_1}$ and $\overrightarrow{v_2}$ are variational symmetries
since they satisfy the condition $~pr^{(1)}(\overrightarrow{v})L+L\mbox{div}(\xi )=0.$

\section*{Acknowledgements}
We would like to thank the referees very much for
their valuable comments and suggestions.

\end{document}